\documentclass[12pt]{article}%
\usepackage{amsmath}
\usepackage{amsfonts}%
\usepackage{amssymb}
\newtheorem{theorem}{Theorem}[section]

\newtheorem{corollary}[theorem]{Corollary}

\newtheorem{definition}[theorem]{Definition}
\newtheorem{example}[theorem]{Example}

\newtheorem{lemma}[theorem]{Lemma}

\newtheorem{proposition}[theorem]{Proposition}
\newtheorem{remark}[theorem]{Remark}

\newcommand{\norm}[1]{\Vert #1 \Vert}

\begin{document}

\title{Hardy algebras associated with $W^{*}$-correspondences (point evaluation and
Schur class functions)}
\author{Paul S. Muhly\thanks{Supported by grants from the U.S. National Science
Foundation and from the U.S.-Israel Binational Science Foundation.}\\Department of Mathematics\\University of Iowa\\Iowa City, IA 52242\\\texttt{muhly@math.uiowa.edu}
\and Baruch Solel\thanks{Supported by the U.S.-Israel Binational Science Foundation
and by the Fund for the Promotion of Research at the Technion}\\Department of Mathematics\\Technion\\32000 Haifa\\Israel\\\texttt{mabaruch@techunix.technion.ac.il}}
\maketitle

\section{Introduction}

This is primarily an exposition of our work in \cite{MSNP} and \cite{MSSchur}
which builds on the theory of tensor algebras over $C^{\ast}$-correspondences
that we developed in \cite{MS98}. Operator tensor algebras (and their
$w^{\ast}$-analogues, which we call \emph{Hardy algebras}) form a rich class
of non-self-adjoint operator algebras that contains a large variety of
operator algebras that have received substantial attention in the literature
in recent years.

Among these algebras are the classical disc algebra $A(\mathbb{D})$, and its
weak closure, $H^{\infty}(\mathbb{T})$; Popescu's non commutative disc
algebras \cite{Po91}, and their weak closures, the free semigroup algebras
studied by Popescu \cite{Po91} and Davidson and Pitts \cite{DP}; quiver
algebras studied by us in \cite{MS99} and by Kribs and Power in \cite{KP};
certain nest algebras; analytic crossed products, studied by Peters \cite{Pe}
and by McAsey and Muhly in \cite{MM}; and others. (We will describe the
construction of tensor and Hardy algebras and give many examples in Section
2.) The theory gives a common approach to the analysis of all these algebras
and has its roots deeply embedded in the model theory of contraction operators
on the one hand and in classical ring theory on the other.

In fact, the theory of contraction operators may be viewed as the theory of
contractive representations of the disc algebra. The representation theory of
the tensor algebras is a natural generalization of this theory that preserves
many of its features. The disc algebra may be viewed as an analytic
generalization of the polynomial algebra in one variable. The interplay
between function theory and the representation theory of the polynomial
algebra has been one of the guiding beacons in model theory for decades
\cite{hB88}. The tensor algebras we analyze are operator algebraic versions of
algebras that generalize polynomial algebras and have been of key importance
in ring theory since 1947 \cite{gH47} and, since 1972, have been a major focus
of attention for analyzing finite dimensional algebras (see \cite{pG72} and
\cite{GR91}). \ (See \cite{pM97} for an extended discussion of the connection
between operator tensor algebras and the theory of finite dimensional algebras.)

Recall that the disc algebra $A(\mathbb{D})$ may be realized as the algebra of
all analytic Toeplitz operators on $l^{2}(\mathbb{N})$ (or on $H^{2}%
(\mathbb{T})$). Popescu generalizes $A(\mathbb{D})$ by considering algebras of
operators on the full Fock space over a Hilbert space $H$ of some dimension,
$n$ say. Let $\mathcal{F}(H)=\mathbb{C}\oplus H\oplus H^{\otimes2}\oplus
\cdots$ denote this Fock space. Then his noncommutative disc algebra of index
$n$ is the norm closed algebra generated by the (left) creation operators.
That is, his algebras are generated by the identity operator and
 operators of
the form $\lambda(\xi)\eta:=\xi\otimes\eta$, where $\xi\in H$ and $\eta
\in\mathcal{F}(H)$. The Fock space may also be written as $l^{2}%
(\mathbb{F}_{+}^{n})$ where $\mathbb{F}_{+}^{n}$ is the free
semigroup on $n$ generators. In this realization, $H$ may be
identified with all the functions supported on the words of length
one and for such a function $\xi$, $\lambda(\xi)$ is just
convolution by $\xi$ on $l^{2}(\mathbb{F}_{+}^{n})$. Observe that
when $n$, the dimension of $H$, is one, then one recovers the disc
algebra $A(\mathbb{D)}$ represented by analytic Toeplitz operators
on $l^{2}(\mathbb{N})$.

 To construct more general tensor algebras
one replaces the Hilbert space $H$ by a correspondence over some
$C^{\ast}$-algebra (or von Neumann algebra) $M$. Roughly, a
correspondence is a bimodule over $M$ that is also equipped with
an $M$-valued inner product. (For a precise definition see Section
2). When $M=\mathbb{C}$ a correspondence over $M$ is just a
Hilbert space.

The tensor algebra associated with $E$, $\mathcal{T}_{+}(E)$, is generated by
creation operators on the Fock space $\mathcal{F}(E)=M\oplus E\oplus
E^{\otimes2}\cdots$ together with a copy of $M$ (formed by diagonal operators
of multiplication, $\varphi_{\infty}(a)$, $a\in M$). It follows from the
results of \cite{MS98} that (completely contractive) representations of
$\mathcal{T}_{+}(E)$ are given by pairs $(T,\sigma)$ where $T:E\rightarrow
B(H)$ is a completely contractive map and $\sigma:M\rightarrow B(H)$ is a
$C^{\ast}$-representation of $M$ that satisfy $T(a\cdot\xi\cdot b)=\sigma
(a)T(\xi)\sigma(b)$ for $a,b\in M$ and $\xi\in E$. (Note that we shall
sometimes use $\varphi$ for the left multiplication on $E$; that is,
$a\cdot\xi$ may be written $\varphi(a)\xi$). Such pairs, $(T,\sigma)$, are
called \emph{covariant representations} of the correspondence $E$. Given
$(T,\sigma)$, one may form the Hilbert space $E\otimes_{\sigma}H$ (see the
discussion following Definition~\ref{hilbertmod}). For $a\in M$,
$\varphi(a)\otimes I$ then defines a bounded operator on this space. The
``complete contractivity'' of $T$ is equivalent to the assertion that the
linear map $\tilde{T}$ defined initially on the balanced algebraic tensor
product $E\otimes H$ by the formula $\tilde{T}(\xi\otimes h):=T(\xi)h$ extends
to an operator of norm at most $1$ on the completion $E\otimes_{\sigma}H$. The
bimodule property of $T$, then, is equivalent to the equation%
\begin{equation}
\tilde{T}(\varphi(a)\otimes I)=\sigma(a)\tilde{T},\label{intertwine0}%
\end{equation}
for all $a\in M$, which means that $\tilde{T}$ intertwines $\sigma$ and the
natural representation of $M$ on $E\otimes_{\sigma}H$ - the composition of
$\varphi$ with Rieffel's induced representation of $\mathcal{L}(E)$ determined
by $\sigma$.

Thus we see that, once $\sigma$ is fixed, the representations of
$\mathcal{T}_{+}(E)$ are parameterized by the elements in the closed unit ball
of the intertwining space $\{\eta\in B(E\otimes_{\sigma}H,H)\mid\eta
(\varphi(\cdot)\otimes I)=\sigma\eta$ and $\left\|  \eta\right\|  \leq1\}$.
Reflecting on this leads one ineluctably to the functional analyst's
imperative: \emph{To understand an algebra, view it as an algebra of functions
on its space of representations. }In our setting, then, we want to think about
$\mathcal{T}_{+}(E)$ as a space of functions on this ball. For reasons that
will be revealed in a minute, we prefer to focus on the adjoints of the
elements in this space. Thus we let $E^{\sigma}=\{\eta\in B(H,E\otimes
_{\sigma}H)\mid\eta\sigma=(\varphi(\cdot)\otimes I)\eta\}$ and we write
$\mathbb{D}((E^{\sigma})^{\ast})$ for the set $\{\eta\in B(E\otimes_{\sigma
}H,H)\mid\eta^{\ast}\in E^{\sigma}$, and $\left\|  \eta\right\|  <1\}$. That
is, $\mathbb{D}((E^{\sigma})^{\ast})$ is the norm-interior of the
representation space consisting of those $(T,\sigma)$ that are ``anchored by
$\sigma$''. One of our interests, then, is to understand the kind of functions
that elements $X$ of $\mathcal{T}_{+}(E)$ determine on $\mathbb{D}((E^{\sigma
})^{\ast})$ via the formula
\[
X(\eta^{\ast})=\sigma\times\eta^{*}(X),
\]
where $\sigma\times\eta^{*}$ is the representation of $\mathcal{T}_{+}(E)$
that is determined by the pair $(\sigma,T)$ with $\tilde{T}=\eta^{*}$.

In the special case when $A=E=\mathbb{C}$ and $\sigma$ is the one-dimensional
representation of $A$ on $\mathbb{C}$, $E^{\sigma}$ is also one-dimensional,
so $\mathbb{D}((E^{\sigma})^{\ast})$ is just the open unit disc in the complex
plane and, for $X\in\mathcal{T}_{+}(E)$, $X(\eta^{\ast})$ is the ordinary
value of $X$ at the complex number $\bar{\eta}$. On the other hand, if
$A=E=\mathbb{C}$, and if $\sigma$ is scalar multiplication on a Hilbert space
$H$ (the only possible representation of $\mathbb{C}$ on $H$), then
$\mathbb{D}((E^{\sigma})^{\ast})$ is the space of \emph{strict }contraction
operators on $H$ and for $\eta^{\ast}\in\overline{\mathbb{D}((E^{\sigma
})^{\ast})}^{\left\|  \cdot\right\|  }$ and $X\in\mathcal{T}_{+}%
(E)=A(\mathbb{D})$, $X(\eta^{\ast})$ is simply the value of $X$ at $\eta
^{\ast}$ defined through the Sz.-Nagy--Foia\c{s} functional calculus
\cite{NF66}. For another example, if $A=\mathbb{C}$, but $E=\mathbb{C}^{n}$,
and if $\sigma$ is scalar multiplication on a Hilbert space $H$, then
$\mathbb{D}((E^{\sigma})^{\ast})$ is the space of row contractions on $H$,
$(T_{1},T_{2},\cdots,T_{n})$, of norm less than $1$; i.e. $\sum T_{i}^{\ast
}T_{i}\leq rI_{H}$ for some $r<1$. In this case, $X(\eta^{\ast})$ is given by
Popescu's functional calculus \cite{Po95}.

In addition to parametrizing certain representations of $\mathcal{T}_{+}(E)$,
$E^{\sigma}$ has another fundamental property: It is itself a $C^{\ast}%
$-correspondence - over \emph{the von Neumann algebra} $\sigma(A)^{\prime}$.
Indeed, it is not difficult to see that $E^{\sigma}$ becomes a bimodule over
$\sigma(A)^{\prime}$ via the formulae: $a\cdot\eta=(I_{E}\otimes a)\eta$ and
$\eta\cdot a=\eta a$, $\eta\in E^{\sigma}$, $a\in\sigma(A)^{\prime}$. Further,
if $\eta$ and $\zeta$ are in $E^{\sigma}$, then the product $\eta^{\ast}\zeta$
lies in the commutant $\sigma(A)^{\prime}$ and defines a $\sigma(A)^{\prime}%
$-valued inner product $\langle\eta,\zeta\rangle$ making $E^{\sigma}$ a
$C^{\ast}$-correspondence. In fact, since $E^{\sigma}$ is a weakly closed
space of operators, it has certain topological properties making it what we
call a $W^{\ast}$-correspondence \cite{MSNP}. It is because $E^{\sigma}$ is a
$W^{\ast}$-correspondence over $\sigma(A)^{\prime}$ that we focus on it, when
studying representations of $\mathcal{T}_{+}(E)$, rather than on its space of
adjoints. While $E^{\sigma}$ plays a fundamental role in our study of quantum
Markov processes \cite{MSQMP}, its importance here - besides providing a space
on which to ``evaluate'' elements of $\mathcal{T}_{+}(E)$ - lies in the fact
that a certain natural representation of $E^{\sigma}$ generates the commutant
of the representation of $\mathcal{T}_{+}(E)$ obtained by ``inducing $\sigma$
up to'' $\mathcal{L}(\mathcal{F}(E))$. (See Theorem~\ref{commutant}).

It is primarily because of this commutant theorem that we cast our work in
this paper entirely in terms of $W^{\ast}$-correspondences. That is, we work
with von Neumann algebras $M$ and $W^{\ast}$-correspondences $E$ over them. We
still form the Fock space $\mathcal{F}(E)$ and the tensor algebra
$\mathcal{T}_{+}(E)$ over $E$, but because $\mathcal{F}(E)$ is a $W^{\ast}%
$-correspondence over $M$, the space $\mathcal{L}(\mathcal{F}(E))$ is a von
Neumann algebra. We call the $w^{\ast}$-closure of $\mathcal{T}_{+}(E)$ in
$\mathcal{L}(\mathcal{F}(E))$ the \emph{Hardy algebra} of $E$ and denote it by
$H^{\infty}(E)$. This is our principal object of study. In the case when
$M=E=\mathbb{C}$, $H^{\infty}(E)$ it the classical $H^{\infty}(\mathbb{T})$
(viewed as analytic Toeplitz operators).

As we will see in Lemma~\ref{contraction}, given a faithful normal
representation $\sigma$ of $M$ on a Hilbert space $H$, we may also evaluate
elements in $H^{\infty}(E)$ at points in $\mathbb{D}((E^{\sigma})^{\ast})$
(since the representation associated with a point in the open unit ball
extends from $\mathcal{T}_{+}(E)$ to $H^{\infty}(E)$). That is, elements in
$H^{\infty}(E)$ may be viewed as functions on $\mathbb{D}((E^{\sigma})^{\ast
})$, also. Further, when $H^{\infty}(E)$ is so represented, one may study the
``value distribution theory'' of these functions. In this context, we
establish two capstone results from function theory: The first, \cite[Theorem
5.3]{MSNP} is presented as Theorem~\ref{NP} below, generalizes the
Nevanlinna-Pick interpolation theorem. It asserts that given two $k$-tuples of
operators in $B(H)$ (where $H$ is the representation space of $\sigma$),
$B_{1},B_{2},\cdots,B_{k}$, and $C_{1},C_{2},\cdots,C_{k}$, and given points
$\eta_{1},\eta_{2},\cdots,\eta_{k}$ in $\mathbb{D}((E^{\sigma})^{\ast})$, one
may find an element $X$ in $H^{\infty}(E)$ of norm at most one such that
\[
B_{i}X(\eta_{i}^{\ast})=C_{i},
\]
for all $i$ if and only if a certain matrix of \emph{maps}, which resembles
the classical Pick matrix, represents a completely positive operator. This
result captures numerous theorems in the literature that go under the name of
generalized Nevanlinna-Pick theorems. Our proof of the theorem (in
\cite{MSNP}) uses a commutant lifting theorem that we proved in \cite{MS98}.
In the context of model theory, it was Sarason who introduced the use of
commutant lifting to prove the interpolation theorem (\cite{S67}). More
recently, a number of authors have been studying interpolation problems in the
context of reproducing kernel Hilbert spaces. (See \cite{BTV} and \cite{AM}).

Our second capstone result is a generalization of Schwartz's lemma (see
Theorem~\ref{Theorem3.5}). It follows from our Nevanlinna-Pick theorem that an
element $X$ in $H^{\infty}(E)$ of norm at most one defines a ``Pick-type''
matrix of maps that represents a completely positive map. In fact, the matrix
is defined using the values of $X$ on $\mathbb{D}((E^{\sigma})^{\ast})$. Given
an arbitrary operator-valued function $Z$ on $\mathbb{D}((E^{\sigma})^{\ast}%
)$, one may define a matrix of maps in a similar way. We say that $Z$ is a
\emph{Schur class operator function} if this matrix defines a completely
positive map. (See Definition~\ref{schur} for a precise statement).
Theorem~\ref{NP} then shows that the function $\eta^{\ast}\mapsto X(\eta
^{\ast})$ is a Schur class operator function for $X$ in the closed unit ball
of $H^{\infty}(E)$. In fact, we show in Theorem~\ref{schurealiz} that every
Schur class operator function arises in this way and that every such function
(with values in, say, $B(\mathcal{E})$) may be represented in the form
$Z(\eta^{\ast})=A+B(I-L_{\eta}^{\ast}D)^{-1}L_{\eta}^{\ast}C$ where $A,B,C$
and $D$ are the entries of a $2\times2$ block matrix representing a
coisometric operator $V$ from $\mathcal{E}\oplus H$ to $\mathcal{E}%
\oplus(E^{\sigma}\otimes H)$ (for some auxiliary Hilbert space $H$) with a
certain intertwining property and $L_{\eta}$ is the operator from $H$ to
$E^{\sigma}\otimes H$ that maps $h$ to $\eta\otimes h$. Borrowing terminology
from the classical function theory on the unit disc $\mathbb{D}$, we call such
a representation a \emph{realization} of $Z$ and we call the coisometry $V$ a
\emph{colligation}. (In general, $V$ is a coisometry but, under a mild
assumption, it may be chosen to be unitary.)

These results, together with our work on canonical models in \cite{MSCM},
represent a generalization of some of the essential ingredients of a program
that has been developed successfully in model theory - the interaction between
operator theory and function theory on the unit disc $\mathbb{D}$ - and has
been generalized in various ways to the polydisc and the ball in
$\mathbb{C}^{n}$. This program sets up (essentially) bijective correspondences
connecting the theory of unitary colligations (and their unitary systems), the
Sz-Nagy-Foias functional model theory for contraction operators and the
discrete-time Lax-Phillips scattering theory. Each theory produces a
contractive operator-valued function (called the transfer function of the
system, the characteristic operator function of the completely non unitary
contraction or the scattering function for the scattering system) from which
one can recover the original object (the system or the contraction) up to
unitary equivalence. For more details, see the works of Ball (\cite{B}), Ball
and Vinnikov (\cite{BV}), Ball, Trent and Vinnikov (\cite{BTV}) and the
references there.

We shall not discuss the program in detail here but we note that
Theorem~\ref{schurealiz} below is the generalization, to our context, of
Theorem 2.1 of \cite{B} or Theorem 2.1 of \cite{BTV}. Here the elements of
$H^{\infty}(E)$ play the role of multipliers and the disc $\mathbb{D}$ in
$\mathbb{C}$ is replaced by the open unit ball of $(E^{\sigma})^{\ast}$.

We also note that the canonical models for contraction operators are replaced,
in our setting, by canonical models for representations of $H^{\infty}(E)$.
This theory was developed in \cite{MSCM} for completely noncoisometric
representations (generalizing results of Popescu in \cite{Po89}) and it is
shown there that the characteristic operator function for such a
representation has a realization associated with a unitary colligation.

In the next section we set the stage by defining our basic constructions,
presenting examples and emphasizing the roles of duality and point evaluation
in the theory.

Section 3 deals with the Nevanlinna-Pick theorem and Section 4 with Schur
class operator functions.

\section{Preliminaries: $W^{\ast}$-correspondences and Hardy algebras}

We shall follow Lance \cite{L94} for the general theory of Hilbert $C^{\ast}%
$-modules that we use. Let $A$ be a $C^{\ast}$-algebra and let $E$ be a right
module over $A$ endowed with a bi-additive map $\langle\cdot,\cdot
\rangle:E\times E\rightarrow A$ (referred to as an $A$-valued inner product)
such that, for $\xi,\eta\in E$ and $a\in A$, $\langle\xi,\eta a\rangle
=\langle\xi,\eta\rangle a$, $\langle\xi,\eta\rangle^{\ast}=\langle\eta
,\xi\rangle$, and $\langle\xi,\xi\rangle\geq0$, with $\langle\xi,\xi\rangle=0$
only when $\xi=0$. If $E$ is complete in the norm $\norm{\xi}
:=\norm{\langle \xi,\xi \rangle} ^{1/2}$, the $E$ is called a (right)
\emph{Hilbert }$C^{\ast}$\emph{-module} over $A$. We write $\mathcal{L}(E)$
for the space of continuous, adjointable, $A$-module maps on $E$; that is
every element of $\mathcal{L}(E)$ is continuous and if $X\in\mathcal{L}(E)$,
then there is an element $X^{\ast}\in\mathcal{L}(E)$ that satisfies $\langle
X^{\ast}\xi,\eta\rangle=\langle\xi,X\eta\rangle$. The element $X^{\ast}$ is
unique and $\mathcal{L}(E)$ is a $C^{\ast}$-algebra with respect to the
involution $X\rightarrow X^{\ast}$ and the operator norm. If $M$ is a von
Neumann algebra and if $E$ is a Hilbert $C^{\ast}$-module over $M$, then $E$
is said to be \emph{self dual} in case every continuous $M$-module map from
$E$ to $M$ is given by an inner product with an element of $E$. If $E$ is a
self dual Hilbert module over $M$, then $\mathcal{L}(E)$ is a $W^{\ast}%
$-algebra and coincides with \emph{all }the bounded linear maps on $E$
\cite{wP73}.

A \emph{$C^{\ast}$-correspondence} over a $C^{\ast}$-algebra $A$ is a Hilbert
$C^{\ast}$-module $E$ over $A$ endowed with a structure of a left module over
$A$ via a *-homomorphism $\varphi:A\rightarrow\mathcal{L}(E)$. When dealing
with a specific $C^{\ast}$-correspondence $E$ over a $C^{\ast} $-algebra $A$,
it will be convenient to suppress the $\varphi$ in formulas involving the left
action and simply write $a\xi$ or $a\cdot\xi$ for $\varphi(a)\xi$. \ This
should cause no confusion in context.

Having defined a left action on $E$, we are allowed to form balanced tensor
products. Given two correspondences $E$ and $F$ over the $C^{\ast}$-algebra
$A$ one may define an $A$-valued inner product on the balanced tensor product
$E\otimes_{A}F$ by the formula
\[
\langle\xi_{1}\otimes\eta_{1},\xi_{2}\otimes\eta_{2}\rangle_{E\otimes_{B}%
F}:=\langle\eta_{1},\varphi(\langle\xi_{1},\xi_{2}\rangle_{E})\eta_{2}%
\rangle_{F}.
\]
The Hausdorff completion of this bimodule is again denoted by $E\otimes_{A}F $
and is called the \emph{tensor product} of $E$ and $F$.

\begin{definition}
\label{hilbertmod} Let $M$ be a von Neumann algebra and let $E$ be a Hilbert
$C^{\ast}$-module over $M$. Then $E$ is called a \emph{Hilbert }$W^{\ast}%
$\emph{-module} over $M$ in case $E$ is self-dual. The module $E$ is called a
$W^{\ast}$\emph{-correspondence over }$M$ in case $E$ is a self-dual $C^{\ast}
$-correspondence over $M$ such that the $\ast$-homomorphism $\varphi: M
\rightarrow\mathcal{L}(E)$ giving the left module structure on $E$ is normal.
\end{definition}

It is evident that the tensor product of two $W^{\ast}$-correspondences is
again a $W^{\ast}$-correspondence. Note also that, given a $W^{\ast}%
$-correspondence $E$ over $M$ and a Hilbert space $H$ equipped with a normal
representation $\sigma$ of $M$, we may form the Hilbert space $E\otimes
_{\sigma}H$ (by defining $\langle\xi_{1}\otimes h_{1},\xi_{2}\otimes
h_{2}\rangle=\langle h_{1},\sigma(\langle\xi_{1},\xi_{2}\rangle)h_{2}\rangle
$). Then, given an operator $X\in\mathcal{L}(E)$ and an operator $S\in
\sigma(N)^{\prime}$, the map $\xi\otimes h\mapsto X\xi\otimes Sh$ defines a
bounded operator on $E\otimes_{\sigma}H$ denoted by $X\otimes S$. When $S=I$
and $X=\varphi(a)$, $a\in M$, we get a representation of $M$ on this space.

Observe that if $E$ is a $W^{\ast}$-correspondence over a von Neumann algebra
$M$, then each of the tensor powers of $E$, viewed as a $C^{\ast}%
$-correspondence over $M$ in the usual way, is in fact a $W^{\ast}%
$-correspondence over $M$ and so, too, is the full Fock space $\mathcal{F}%
(E)$, which is defined to be the direct sum $M\oplus E\oplus E^{\otimes
2}\oplus\cdots$, with its obvious structure as a right Hilbert module over $M
$ and left action given by the map $\varphi_{\infty}$, defined by the formula
$\varphi_{\infty}(a):=diag\{a,\varphi(a),\varphi^{(2)}(a),\varphi
^{(3)}(a),\cdots\}$, where for all $n$, $\varphi^{(n)}(a)(\xi_{1}\otimes
\xi_{2}\otimes\cdots\xi_{n})=(\varphi(a)\xi_{1})\otimes\xi_{2}\otimes\cdots
\xi_{n}$, $\xi_{1}\otimes\xi_{2}\otimes\cdots\xi_{n}\in E^{\otimes n}$. The
\emph{tensor algebra} over $E$, denoted $\mathcal{T}_{+}(E)$, is defined to be
the norm-closed subalgebra of $\mathcal{L}(\mathcal{F}(E))$ generated by
$\varphi_{\infty}(M)$ and the \emph{creation operators} $T_{\xi}$, $\xi\in E$,
defined by the formula $T_{\xi}\eta=\xi\otimes\eta$, $\eta\in\mathcal{F}(E) $.
We refer the reader to \cite{MS98} for the basic facts about $\mathcal{T}%
_{+}(E)$.

\begin{definition}
\label{Hinfty}Given a $W^{\ast}$-correspondence $E$ over the von Neumann
algebra $M$, the ultraweak closure of the tensor algebra of $E$,
$\mathcal{T}_{+}(E)$, in the $w^{\ast}$-algebra $\mathcal{L}(\mathcal{F}(E))$,
will be called the \emph{Hardy Algebra of }$E$, and will be denoted
$H^{\infty}(E)$.
\end{definition}

\begin{example}
\label{ex0.1} If $M=E=\mathbb{C}$ then $\mathcal{F}(E)$ may be identified with
$H^{2}(\mathbb{T})$. The tensor algebra in this setting is isomorphic to the
disc algebra $A(\mathbb{D})$ and the Hardy algebra is the classical Hardy
algebra $H^{\infty}(\mathbb{T})$.
\end{example}

\begin{example}
\label{ex0.2} If $M=\mathbb{C}$ and $E=\mathbb{C}^{n}$, then $\mathcal{F}(E)$
may be identified with the space $l_{2}(\mathbb{F}_{n}^{+})$ where
$\mathbb{F}_{n}^{+}$ is the free semigroup on $n$ generators. The tensor
algebra then is what Popescu refers to as the ``non commutative disc algebra''
$\mathcal{A}_{n}$ and the Hardy algebra is its $w^{\ast}$-closure. It was
studied by Popescu (\cite{Po01}) and by Davidson and Pitts who denoted it by
$\mathcal{L}_{n}$ (\cite{DP}).
\end{example}

\begin{example}
\label{ex0.3} Let $M$ be a von Neumann algebra and let $\alpha$ be a unital,
injective, normal $^{\ast}$-endomorphism on $M$. The correspondence $E$
associated with $\alpha$ is equal to $M$ as a vector space. The right action
is by multiplication, the $M$-valued inner product is $\langle a,b\rangle
=a^{\ast}b$ and the left action is given by $\alpha$; i.e. $\varphi
(a)b=\alpha(a)b$. We write $_{\alpha}M$ for $E$. It is easy to check that
$E^{\otimes n}$ is isomorphic to $_{\alpha^{n}}M$. The Hardy algebra in this
case is called the \emph{non-selfadjoint crossed product} of $M$ by $\alpha$
and will be written $M\rtimes_{\alpha}\mathbb{Z}_{+}$. This algebra is also
called an \emph{analytic crossed product}, at least when $\alpha$ is an
automorphism. It is related to the algebras studied in \cite{MM} and
\cite{Pe}. If we write $w$ for $T_{1}$ (where $1$ is the identity of $M$
viewed as an element of $E$), then the algebra is generated by $w$ and
$\varphi_{\infty}(M)$ and every element $X$ in the algebra has a formal
``Fourier'' expression
\[
X=\sum_{n=0}w^{n}b_{n}%
\]
where $b_{n}\in\varphi_{\infty}(M)$. This Fourier expression is actually
Ceasaro-summable to $X$ in the ultraweak topology on $H^{\infty}(E)$
\cite{MSNP}, but we don't need these details in the present discussion.
\end{example}

\begin{example}
\label{Upper} Here we set $M$ to be the algebra $l_{\infty}(\mathbb{Z})$ and
let $\alpha$ be the automorphism defined by $(\alpha(g))_{i}=g_{i-1}$. Write
$E $ for the correspondence $_{\alpha}M$ as in Example~\ref{ex0.3}. Another,
isomorphic, way to describe $E$ is to let $M$ be the algebra $\mathcal{D}$ of
all diagonal operators on $l^{2}(\mathbb{Z})$ , let $U$ be the shift defined
by $Ue_{k}=e_{k-1}$ (where $\{e_{k}\}$ is the standard basis), and set
$E=U\mathcal{D} \subseteq B(l^{2}(\mathbb{Z}))$. The left and right actions on
$E$ are defined simply by operator multiplications and the inner product is
$\langle UD_{1},UD_{2} \rangle=D_{1}^{*}D_{2}$. It is easy to check that these
correspondences are indeed isomorphic and the Hardy algebra $H^{\infty}(E)$ is
(completely isometrically isomorphic to) the algebra $\mathcal{U}$ of all
operators in $B(l_{2}(\mathbb{Z}))$ whose matrix (with respect to the standard
basis) is upper triangular.
\end{example}

\begin{example}
\label{cp} Suppose that $\Theta$ is a normal, contractive, completely positive
map on a von Neumann algebra $M$. Then we may associate with it the
correspondence $M\otimes_{\Theta}M$ obtained by defining the $M$-valued inner
product on the algebraic tensor product $M\otimes M$ via the formula $\langle
a\otimes b,c\otimes d\rangle=b^{\ast}\theta(a^{\ast}c)d$ and completing. (The
bimodule structure is by left and right multiplications). This correspondence
was used by Popa (\cite{P86}), Mingo (\cite{jM}), Anantharam-Delarouche
(\cite{AD90}) and others to study the map $\Theta$. If $\Theta$ is an
endomorphism this correspondence is the one described in example~\ref{ex0.3}.
\end{example}

\begin{example}
\label{ex0.4} Let $M$ be $D_{n}$, the diagonal $n\times n$ matrices and $E$ be
the set of all $n\times n$ matrices $A=(a_{ij})$ with $a_{ij}=0$ unless
$j=i+1$ with the inner product $\langle A,B\rangle=A^{\ast}B$ and the left and
right actions given by matrix multiplication. Then the Hardy algebra is
isomorphic to $T_{n}$, the $n\times n$ upper triangular matrices. In fact, a
similar argument works to show that, for every \emph{finite} nest of
projections $\mathcal{N}$ on a Hilbert space $H$, the nest algebra
$alg\mathcal{N}$ (i.e. the set of all operators on $H$ leaving the ranges of
the projections in $\mathcal{N}$ invariant ) may be written as $H^{\infty}(E)$
for some $W^{\ast}$-correspondence $E$.
\end{example}

\begin{example}
\label{quiver} (Quiver algebras) Let $\mathcal{Q}$ be a directed graph on the
set $V$ of vertices. For simplicity we assume that both $V$ and $\mathcal{Q}$
are finite sets and view each $\alpha\in\mathcal{Q}$ as an arrow from
$s(\alpha)$ (in $V$) to $r(\alpha)$ (in $V$). Let $M$ be $C(V)$ (a von Neumann
algebra) and $E$ (or $E(\mathcal{Q})$) be $C(\mathcal{Q})$. Define the
$M$-bimodule structure on $E$ as follows: for $f\in E$, $\psi\in M $ and
$\alpha\in\mathcal{Q}$,
\[
(f\psi)(\alpha)=f(\alpha)\psi(s(\alpha)) ,
\]
and
\[
(\psi f)(\alpha)=\psi(r(\alpha))f(\alpha).
\]
The $M$-valued inner product is given by the formula
\[
\langle f,g \rangle(v)=\sum_{s(\alpha)=v} \overline{f(\alpha)}g(\alpha) ,
\]
for $f,g \in E$ and $v\in V$. The algebra $H^{\infty}(E)$ in this case will be
written $H^{\infty}(\mathcal{Q})$ and is the $\sigma$-weak closure of
$\mathcal{T}_{+}(E(\mathcal{Q}))$. Viewing both algebras as acting on the Fock
space, one sees that they are generated by a set $\{S_{\alpha} : \alpha
\in\mathcal{Q} \}$ of partial isometries (here $S_{\alpha}=T_{\delta_{\alpha}%
}$ where $\delta_{\alpha}$ is the function in $C(\mathcal{Q})$ which is $1$ at
$\alpha$ and $0$ otherwise) and a set $\{P_{v} : v\in V\}$ of projections
(i.e. the generators of $\varphi_{\infty}(M)$) satisfying the following conditions.

\begin{enumerate}
\item[(i)] $P_{v}P_{u}=0 $ if $u\neq v$,

\item[(ii)] $S_{\alpha}^{*}S_{\beta}=0$ if $\alpha\neq\beta$

\item[(iii)] $S_{\alpha}^{*}S_{\alpha}=P_{s(\alpha)}$ and

\item[(iv)] $\sum_{r(\alpha)=v} S_{\alpha}S_{\alpha}^{*} \leq P_{v}$ for all
$v\in V$.
\end{enumerate}

These algebras were studied in \cite{pM97} and \cite{MS99}, and also in
\cite{KP} where they were called free semigroupoid algebras.
\end{example}

\subsection{Representations}

In most respects, the representation theory of $H^{\infty}(E)$ follows the
lines of the representation theory of $\mathcal{T}_{+}(E)$. However, there are
some differences that will be important to discuss here. To help illuminate
these, we need to review some of the basic ideas from \cite{MS98, MS99, MSNP}.

A representation $\rho$ of $H^{\infty}(E)$ (or of $\mathcal{T}_{+}(E)$) on a
Hilbert space $H$ is completely determined by what it does to the generators.
Thus, from a representation $\rho$ we obtain two maps: a map $T:E\rightarrow
B(H)$, defined by $T(\xi)=\rho(T_{\xi})$, and a map $\sigma:M\rightarrow
B(H)$, defined by $\sigma(a)=\rho(\varphi_{\infty}(a))$. Analyzing the
properties of $T$ and $\sigma$ one is lead to the following definition.

\begin{definition}
\label{Definition1.12}Let $E$ be a $W^{\ast}$-correspondence over a von
Neumann algebra $M$. Then a \emph{completely contractive covariant
representation }of $E$ on a Hilbert space $H$ is a pair $(T,\sigma)$, where

\begin{enumerate}
\item $\sigma$ is a normal $\ast$-representation of $M$ in $B(H)$.

\item $T$ is a linear, completely contractive map from $E$ to $B(H)$ that is
continuous in the $\sigma$-topology of \cite{BDH88} on $E$ and the ultraweak
topology on $B(H).$

\item $T$ is a bimodule map in the sense that $T(S\xi R)=\sigma(S)T(\xi
)\sigma(R)$, $\xi\in E$, and $S,R\in M$.
\end{enumerate}
\end{definition}

It should be noted that there is a natural way to view $E$ as an operator
space (by viewing it as a subspace of its linking algebra) and this defines
the operator space structure of $E$ to which the Definition
\ref{Definition1.12} refers when it is asserted that $T$ is completely contractive.

As we noted in the introduction and developed in \cite[Lemmas 3.4--3.6]{MS98}
and in \cite{MSNP}, if a completely contractive covariant representation,
$(T,\sigma)$, of $E$ in $B(H)$ is given, then it determines a contraction
$\tilde{T}:E\otimes_{\sigma}H\rightarrow H$ defined by the formula $\tilde
{T}(\eta\otimes h):=T(\eta)h$, $\eta\otimes h\in E\otimes_{\sigma}H$. The
operator $\tilde{T}$ satisfies
\begin{equation}
\tilde{T}(\varphi(\cdot)\otimes I)=\sigma(\cdot)\tilde{T}.\label{covariance}%
\end{equation}
In fact we have the following lemma from \cite[Lemma 2.16]{MSNP}.

\begin{lemma}
\label{CovRep}The map $(T,\sigma)\rightarrow\tilde{T}$ is a bijection between
all completely contractive covariant representations $(T,\sigma)$ of $E$ on
the Hilbert space $H$ and contractive operators $\tilde{T}:E\otimes_{\sigma
}H\rightarrow H$ that satisfy equation (\ref{covariance}). Given $\sigma$ and
a contraction $\tilde{T}$ satisfying the covariance condition
(\ref{covariance}), we get a the completely contractive covariant
representation $(T,\sigma)$ of $E$ on $H$ by setting $T(\xi)h:=\tilde{T}%
(\xi\otimes h)$.
\end{lemma}

The following theorem shows that every completely contractive representation
of the tensor algebra $\mathcal{T}_{+}(E)$ is given by a pair $(T,\sigma)$ as
above or, equivalently, by a contraction $\tilde{T}$ satisfying
(\ref{covariance}).

\begin{theorem}
\label{representation} (\cite[Theorem 3.10]{MS98}) Let $E$ be a $W^{\ast}%
$-correspondence over a von Neumann algebra $M$. To every completely
contractive covariant representation, $(T,\sigma)$, of $E$ there is a unique
completely contractive representation $\rho$ of the tensor algebra
$\mathcal{T}_{+}(E)$ that satisfies
\[
\rho(T_{\xi})=T(\xi)\;\;\;\xi\in E
\]
and%
\[
\rho(\varphi_{\infty}(a))=\sigma(a)\;\;\;a\in M.
\]
The map $(T,\sigma)\mapsto\rho$ is a bijection between the set of all
completely contractive covariant representations of $E$ and all completely
contractive (algebra) representations of $\mathcal{T}_{+}(E)$ whose
restrictions to $\varphi_{\infty}(M)$ are continuous with respect to the
ultraweak topology on $\mathcal{L}(\mathcal{F}(E))$.
\end{theorem}

\begin{definition}
\label{integratedform}If $(T,\sigma)$ is a completely contractive covariant
representation of a $W^{\ast}$-correspondence $E$ over a von Neumann algebra
$M$, we call the representation $\rho$ of $\mathcal{T}_{+}(E)$ described in
Theorem \ref{representation} the \emph{integrated form} of $(T,\sigma)$ and
write $\rho=\sigma\times T$.
\end{definition}

\begin{example}
\label{0.2rep} In the context of Example~\ref{ex0.2}, $M=\mathbb{C}$ and
$E=\mathbb{C}^{n}$. Then, a completely contractive covariant representation of
$E$ is simply given by a completely contractive map $T:E\rightarrow B(H)$.
Writing $T_{k}=T(e_{k})$, where $e_{k}$ is the standard basis in
$\mathbb{C}^{n}$, and identifying $\mathbb{C}^{n}\otimes H$ with the direct
sum of $n$ copies of $H$, we may write $\tilde{T}$ as a row $(T_{1}%
,T_{2},\ldots,T_{n})$. The condition that $\norm{\tilde{T}} \leq1$ is the
condition (studied by Popescu \cite{Po91} and Davidson and Pitts \cite{DP})
that $\sum T_{i}T_{i}^{\ast}\leq1$. Hence representations of the
noncommutative disc algebras are given by row contractions.
\end{example}

\begin{example}
\label{quiverrep} Consider the setting of Example~\ref{quiver} and let
$V,\mathcal{Q},M$ and $E$ be as defined there. A (completely contractive
covariant) representation of $E$ is given by a representation $\sigma$ of
$M=C(V)$ on a Hilbert space $H$ and by a contractive map $\tilde{T}%
:E\otimes_{\sigma}H\rightarrow H$ satisfying (\ref{covariance}) above. Write
$\delta_{v}$ for the function in $C(V)$ which is $1$ on $v$ and $0$ elsewhere.
The representation $\sigma$ is given by the projections $Q_{v}=\sigma
(\delta_{v})$ whose sum is $I$. For every $\alpha\in\mathcal{Q}$ write
$\delta_{\alpha}$ for the function (on $E$) which is $1$ at $\alpha$ and $0$
elsewhere. Given $\tilde{T}$ as above, we may define maps $T(\alpha)\in B(H)$
by $T(\alpha)h=\tilde{T}(\delta_{\alpha}\otimes h)$ and it is easy to check
that $\tilde{T}\tilde{T}^{\ast}=\sum_{\alpha}T(\alpha)T(\alpha)^{\ast}$ and
$T(\alpha)=Q_{r(\alpha)}T(\alpha)Q_{s(\alpha)}$. Thus to every (completely
contractive) representation of the quiver algebra $\mathcal{T}_{+}%
(E(\mathcal{Q}))$ we associate a family $\{T(\alpha):\alpha\in\mathcal{Q}\}$
of maps on $H$ that satisfy $\sum_{\alpha}T(\alpha)T(\alpha)^{\ast}\leq I$ and
$T(\alpha)=Q_{r(\alpha)}T(\alpha)Q_{s(\alpha)}$. Conversely, every such family
defines a representation by writing $\tilde{T}(f\otimes h)=\sum f(\alpha
)T(\alpha)h$. Thus, representations are indexed by such families. Note that,
in fact, $(\sigma\times T)(S_{\alpha})=T(\alpha)$ and $(\sigma\times
T)(P_{v})=Q_{v}$ (where $S_{\alpha}$ and $P_{v}$ are as in
Example~\ref{quiver}).
\end{example}

\begin{remark}
\label{keyproblem}One of the principal difficulties one faces in dealing with
$\mathcal{T}_{+}(E)$ and $H^{\infty}(E)$ is to decide when the integrated
form, $\sigma\times T$, of a completely contractive covariant representation
$(T,\sigma)$ extends from $\mathcal{T}_{+}(E)$ to $H^{\infty}(E)$. This
problem arises already in the simplest situation, vis. when $M=\mathbb{C}=E$.
In this setting, $T$ is given by a single contraction operator $T(1)$ on a
Hilbert space, $\mathcal{T}_{+}(E)$ ``is'' the disc algebra and $H^{\infty
}(E)$ ``is'' the space of bounded analytic functions on the disc. The
representation $\sigma\times T$ extends from the disc algebra to $H^{\infty
}(E)$ precisely when there is no singular part to the spectral measure of the
minimal unitary dilation of $T(1)$. We are not aware of a comparable result in
our general context but we have some sufficient conditions. One of them is
given in the following lemma. It is not necessary in general.
\end{remark}

\begin{lemma}
\label{contraction} (\cite{MSNP}) If $\norm{\tilde{T}} <1$ then $\sigma\times
T$ extends to a $\sigma$-weakly continuous representation of $H^{\infty}(E)$.
\end{lemma}

Other sufficient conditions are presented in Section 7 of \cite{MSNP}.

\subsection{Duality and point evaluation}

The following definition is motivated by condition (\ref{covariance}) above.

\begin{definition}
Let $\sigma:M\rightarrow B(H)$be a normal representation of the von Neumann
algebra $M$ on the Hilbert space $H$. Then for a $W^{\ast}$-correspondence $E
$ over $M$, the $\sigma$\emph{-dual} of $E$, denoted $E^{\sigma}$, is defined
to be
\[
\{\eta\in B(H,E\otimes_{\sigma}H)\mid\eta\sigma(a)=(\varphi(a)\otimes
I)\eta,\;a \in M \}.
\]
\end{definition}

As we note in the following proposition, the $\sigma$-dual carries a natural
structure of a $W^{\ast}$-correspondence. The reason to define the $\sigma
$-dual using covariance condition which is the ``adjoint'' of condition
(\ref{covariance}) is to get a \emph{right} $W^{\ast}$-module (instead of a
left $W^{\ast}$-module) over $\sigma(M)^{\prime}$.

\begin{proposition}
With respect to the actions of $\sigma(M)^{\prime}$ and the $\sigma
(M)^{\prime}$-valued inner product defined as follows, $E^{\sigma}$ becomes a
$W^{\ast}$-correspondence over $\sigma(M)^{\prime}$: For $a,b\in
\sigma(M)^{\prime}$, and $\eta\in E^{\sigma}$, $a\cdot\eta\cdot b: =(I\otimes
a)\eta b$, and for $\eta,\zeta\in E^{\sigma}$, $\langle\eta,\zeta
\rangle_{\sigma(M)^{\prime}}:=\eta^{\ast}\zeta$.
\end{proposition}

\begin{example}
\label{ex0.6} If $M=E=\mathbb{C}$, $H$ is arbitrary and $\sigma$ is the
representation of $\mathbb{C}$ on $H$, then $\sigma(M)^{\prime}=B(H)$ and
$E^{\sigma}=B(H)$.
\end{example}

\begin{example}
\label{ex0.7} If $\Theta$ is a contractive, normal, completely positive map on
a von Neumann algebra $M$ and if $E=M\otimes_{\Theta}M$ (see Example~\ref{cp}
) then, for every faithful representation $\sigma$ of $M$ on $H$, the $\sigma
$-dual is the space of all bounded operators mapping $H$ into the Stinespring
space $K$ (associated with $\Theta$ as a map from $M$ to $B(H)$) that
intertwine the representation $\sigma$ (on $H$) and the Stinespring
representation $\pi$ (on $K$). This correspondence was proved very useful in
the study of completely positive maps. (See \cite{MSQMP}, \cite{MSS03} and
\cite{MSC}). If $M=B(H)$ this is a Hilbert space and was studied by Arveson
(\cite{Ar89}). Note also that, if $\Theta$ is an endomorphism, then this dual
correspondence is the space of all operators on $H$ intertwining $\sigma$ and
$\sigma\circ\Theta$.
\end{example}

We now turn to define point evaluation. Note that, given $\sigma$ as above,
the operators in $E^{\sigma}$ whose norm does not exceed $1$ are precisely the
adjoints of the operators of the form $\tilde{T}$ for a covariant pair
$(T,\sigma)$. In particular, every $\eta$ in the \emph{open} unit ball of
$E^{\sigma}$ (written $\mathbb{D}(E^{\sigma})$) gives rise to a covariant pair
$(T,\sigma)$ (with $\eta=\tilde{T}^{\ast}$) such that $\sigma\times T$ is a
representation of $H^{\infty}(E)$. Given $X\in H^{\infty}(E)$ we may apply
$\sigma\times T$ to it. The resulting operator in $B(H)$ will be denoted by
$X(\eta^{\ast})$. That is,
\[
X(\eta^{\ast})=(\sigma\times T)(X)
\]
where $\tilde{T}=\eta^{\ast}$.

In this way, we view every element in the Hardy algebra as a ($B(H)$-valued)
function on $\mathbb{D}((E^{\sigma})^{*})$.

\begin{example}
\label{ex0.12} Suppose $M=E=\mathbb{C}$ and $\sigma$ the representation of
$\mathbb{C}$ on some Hilbert space $H$. Then $H^{\infty}(E)=H^{\infty
}(\mathbb{T})$ and (Example~\ref{ex0.6}) $E^{\sigma}$ is isomorphic to $B(H)$.
If $X\in H^{\infty}(E)=H^{\infty}(\mathbb{T})$, so that we may view $X$ with a
bounded analytic function on the open disc in the plane, then for $S\in
E^{\sigma}=B(H)$, it is not hard to check that $X(S^{\ast})$, as defined
above, is the same as the value provided by the Sz.-Nagy-Foia\c{s} $H^{\infty
}$-functional calculus.
\end{example}

Note that, for a given $\eta\in\mathbb{D}(E^{\sigma})$, the map $X \mapsto
X(\eta^{*})$ is a $\sigma$-weakly continuous homomorphism on the Hardy
algebra. Thus, in order to compute $X(\eta^{*})$, it suffices to know its
values on the generators. This is given in the following (easy to verify) lemma.

\begin{lemma}
\label{pteval1} Let $\sigma$ be a faithful normal representation of $M$ on $H$
and for $\xi\in E$ write $L_{\xi}$ for the map from $H$ to $E\otimes_{\sigma
}H$ defined by $L_{\xi}h=\xi\otimes h$. Then, for $\xi\in E$, $a\in M$ and
$\eta\in\mathbb{D}(E^{\sigma})$,

\begin{enumerate}
\item[(i)] $(T_{\xi})(\eta^{*})=\eta^{*} \circ L_{\xi}$, and

\item[(ii)] $(\varphi_{\infty}(a))(\eta^{*})=\sigma(a) $
\end{enumerate}

\noindent (Recall that $\eta^{\ast}$ is a map from $E\otimes_{\sigma}H$ to $H$.)
\end{lemma}

A formula for computing $X(\eta^{\ast})$, without referring to the generators,
will be presented later (Proposition~\ref{formula}).

\begin{example}
\label{ex0.3eval} In the setting of Example~\ref{ex0.3} we may identify the
Hilbert space $E\otimes_{\sigma}H={}_{\alpha}M\otimes_{\sigma}H$ with $H$ via
the unitary operator mapping $a\otimes h$ (in $_{\alpha}M\otimes_{\sigma}H$)
to $\sigma(a)h$. Using this identification, we may identify $E^{\sigma}$ with
$\{\eta\in B(H):\eta\sigma(a)=\sigma(\alpha(a))\eta,\;a\in M\}$.

Applying the lemma \ref{pteval1}, we obtain $w(\eta^{\ast})=T_{1}(\eta^{\ast
})=\eta^{\ast}\circ L_{1}=\eta^{\ast}$ (viewed now as an operator in $B(H)$).
Thus, if $X=\sum w^{n}b_{n}$ (as a formal series), with $b_{n}=\varphi
_{\infty}(a_{n})$ and $\eta\in\mathbb{D}(E^{\sigma})$, then
\[
X(\eta^{\ast})=\sum(\eta^{\ast})^{n}\sigma(a_{n})
\]
with the sum converging in the norm on $B(H)$. (In a sense, this equation
asserts that Ceasaro summability implies Abel summability even in this
abstract setting.)
\end{example}

\begin{example}
\label{uppereval} Let $\mathcal{D}$, $U$ and $E=U\mathcal{D}$ be as in
Example~\ref{Upper}. Let $\sigma$ be the identity representation of
$\mathcal{D}$ on $H=l^{2}(\mathbb{Z})$. The map $V(UD\otimes_{\sigma}h)=Dh$
(for $D\in\mathcal{D},h\in H$) is a unitary operator from $E\otimes_{\sigma}H$
onto $H$ such that, for every $\eta\in E^{\sigma}$, $V\eta\in U^{\ast
}\mathcal{D}$ and, conversely, for every $D\in\mathcal{D}$, $V^{\ast}U^{\ast
}D^{\ast}$ lies in $E^{\sigma}$. We write $\eta_{D}$ for $V^{\ast}U^{\ast
}D^{\ast}$. Recall that the Hardy algebra is $\mathcal{U}$ (the algebra of all
upper triangular operators on $H$). Given $X\in\mathcal{U}$ we shall write
$X_{n}$ for the $n^{th}$ upper diagonal of $X$. A simple computation shows
that, for $D\in\mathcal{D}$ with $\norm{D} <1$,
\[
X(\eta_{D}^{\ast})=\sum_{n=0}^{\infty}U^{n}(U^{\ast}D)^{n}X_{n}.
\]
Note here that, in \cite{ADD90}, the authors defined point evaluations for
operators $X\in\mathcal{U}$. In their setting one evaluates $X$ on the open
unit ball of $\mathcal{D}$ and the values are also in $\mathcal{D}$. Their
formula (for what in \cite{ADP} is called the right point evaluation) is
\[
X^{\Delta}(D)=\sum_{n=0}^{\infty}U^{n}(U^{\ast}D)^{n}X_{n}U^{\ast n}.
\]
(One can also define a left point evaluation). The apparent similarity of the
two formulas above may be deceiving. Note that both their point evaluation and
ours can be defined also for ``block upper triangular '' operators (acting on
$l^{2}(\mathbb{Z},K)$ for some Hilbert space $K$). But, in that case, the
relation between the two formulas is no longer clear. In fact, our point
evaluation is multiplicative (that is, $(XY)(\eta^{\ast})=X(\eta^{\ast}%
)Y(\eta^{\ast})$) while theirs is not. On the other hand, their point
evaluation is ``designed'' to satisfy the property that, for $X\in\mathcal{U}$
and $D\in\mathcal{D}$, $(X-X^{\Delta}(D))(U-D)^{-1}\in\mathcal{U}$
(\cite[Theorem 3.4]{ADD90}) . For our point evaluation (in the general
setting), it is not even clear how to state such a property.
\end{example}

\begin{example}
\label{quivereval} (Quiver algebras) Let $\mathcal{Q}$ be a quiver as in
Example~\ref{quiver} and write $E(\mathcal{Q})$ for the associated
correspondence. We fix a faithful representation $\sigma$ of $M=C(V)$ on $H$.
As we note in Example~\ref{quiverrep} , this gives a family $\{Q_{v}\}$ of
projections whose sum is $I$ (and, as $\sigma$ is faithful, none is $0$).
Write $H_{v}$ for the range of $Q_{v}$. Then $\sigma(M)^{\prime}=\oplus
_{v}B(H_{v})$ and we write elements there as functions $\psi$ defined on $V$
with $\psi(v)\in B(H_{v})$. To describe the $\sigma$-dual of $E$ we can use
Example 3.4 in \cite{MSNP}. We may also use the description of the maps
$\tilde{T}$ in Example~\ref{quiverrep} because every $\eta$ in the closed unit
ball of $E^{\sigma}$ is $\tilde{T}^{\ast}$ for some representation
$(\sigma,T)$ of $E$. Using this, we may describe an element $\eta$ of
$E^{\sigma}$ as a family of $B(H)$-valued operators $\{\eta(\beta):\beta
\in\mathcal{Q}^{-1}\}$ where $\mathcal{Q}^{-1}$ is the quiver obtained from
$\mathcal{Q}$ by reversing all arrows. The $\sigma(M)^{\prime}$-module
structure of $E^{\sigma}$ is described as follows. For $\eta\in E^{\sigma}$,
$\psi\in\sigma(M)^{\prime}$ and $\beta\in\mathcal{Q}^{-1}$,
\[
(\eta\psi)(\beta)=\eta(\beta)\psi(s(\beta)),
\]
and
\[
(\psi\eta)(\beta)=\psi(r(\beta))\eta(\beta).
\]
The $\sigma(M)^{\prime}$-valued inner product is given by the formula
\[
\langle\eta,\zeta\rangle(v)=\sum_{s(\beta)=v}\eta(\beta)^{\ast}\zeta(\beta),
\]
for $\eta,\zeta\in E^{\sigma}$ and $v\in V$.

Recall that the quiver algebra is generated by a set of partial isometries
$\{S_{\alpha}\}$ and projections $\{P_{v}\}$ (see Example~\ref{quiver}). If
$\sigma$ is given and $\eta^{\ast}=\tilde{T}$ lies in the open unit ball of
$(E^{\sigma})^{\ast}$ and $\tilde{T}$ is given by a row contraction
$(T(\alpha))$ (as in Example~\ref{quiverrep}), then the point evaluation for
the generators is defined by $S_{\alpha}(\eta^{\ast})=T(\alpha)=\eta
(\alpha^{-1})^{\ast}$ and $P_{v}(\eta^{\ast})=Q_{v}$. For a general $X\in
H^{\infty}(\mathcal{Q})$, $X(\eta^{\ast})$ is defined by the linearity,
multiplicativity and $\sigma$-weak continuity of the map $X\mapsto
X(\eta^{\ast})$.
\end{example}

We turn now to some general results concerning the $\sigma$-dual. First, the
term ``dual" that we use is justified by the following result.

\begin{theorem}
\label{dual} (\cite[Theorem 3.6]{MSNP}) Let $E$ be a $W^{\ast}$-correspondence
over $M$ and let $\sigma$ be a faithful, normal representation of $M$ on $H$.
If we write $\iota$ for the identity representation of $\sigma(M)^{\prime}$
(on $H$) then one may form the $\iota$-dual of $E^{\sigma}$ and we have
\[
(E^{\sigma})^{\iota}\cong E.
\]
\end{theorem}

The following lemma summarizes Lemmas 3.7 and 3.8 of \cite{MSNP} and shows
that the operation of taking duals behaves nicely with respect to direct sums
and tensor products.

\begin{lemma}
\label{tensorsum} Given $W^{\ast}$-correspondences $E$,$E_{1}$ and $E_{2}$
over $M$ and a faithful representation $\sigma$ of $M$ on $H$, we have

\begin{enumerate}
\item[(i)] $(E_{1} \oplus E_{2})^{\sigma}\cong E_{1}^{\sigma} \oplus
E_{2}^{\sigma}.$

\item[(ii)] $(E_{1} \otimes E_{2})^{\sigma}\cong E_{2}^{\sigma} \otimes
E_{1}^{\sigma}.$

\item[(iii)] $\mathcal{F}(E)^{\sigma}\cong\mathcal{F}(E^{\sigma}).$

\item[(iv)] The map $\eta\otimes h \mapsto\eta(h)$ induces a unitary operator
from $E^{\sigma}\otimes_{\iota}H$ onto $E\otimes_{\sigma}H$.

\item[(v)] Applying item (iv) above to $\mathcal{F}(E)$ in place of $E$, we
get a unitary operator $U$ from $\mathcal{F}(E^{\sigma})\otimes H$ onto
$\mathcal{F}(E)\otimes H$.
\end{enumerate}
\end{lemma}

Although $H^{\infty}(E)$ was defined as a subalgebra of $\mathcal{L}%
(\mathcal{F}(E))$ it is often useful to consider a (faithful) representation
of it on a Hilbert space. Given a faithful, normal, representation $\sigma$ of
$M$ on $H$ we may ``induce'' it to a representation of the Hardy algebra. To
do this, we form the Hilbert space $\mathcal{F}(E)\otimes_{\sigma}H$ and
write
\[
Ind(\sigma)(X)=X\otimes I,\;\;X\in H^{\infty}(E).
\]
(in fact, this is well defined for every $X$ in $\mathcal{L}(\mathcal{F}(E))$.
Such representations were studied by M. Rieffel in \cite{Rie}). $Ind(\sigma)$
is a faithful representation and is an homeomorphism with respect to the
$\sigma$-weak topologies. Similarly one defines $Ind(\iota)$, a representation
of $H^{\infty}(E^{\sigma})$. The following theorem shows that, roughly
speaking, the algebras $H^{\infty}(E)$ and $H^{\infty}(E^{\sigma})$ are the
commutant of each other.

\begin{theorem}
\label{commutant}\cite[Theorem 3.9]{MSNP} With the operator $U$ as in part (v)
of Lemma~\ref{tensorsum}, we have
\[
U^{\ast}(Ind(\iota)(H^{\infty}(E^{\sigma})))U=(Ind(\sigma)(H^{\infty
}(E)))^{\prime}%
\]
and, consequently,
\[
(Ind(\sigma)(H^{\infty}(E)))^{\prime\prime}=Ind(\sigma)(H^{\infty}(E)).
\]
\end{theorem}

We may now use the notation set above to present a general formula for point
evaluation. For its proof, see \cite[Proposition 5.1]{MSNP}.

\begin{proposition}
\label{formula} If $\sigma$ is a faithful normal representation of $M$ on $H$,
let $\iota_{H}$ denote the imbedding of $H$ into $\mathcal{F}(E^{\sigma
})\otimes H$ and write $P_{k}$ for the projection of $\mathcal{F}(E^{\sigma
})\otimes H$ onto $(E^{\sigma})^{\otimes k}\otimes H$. Also, for $\eta
\in\mathbb{D}(E^{\sigma})$ and $k\geq1$, note that $\eta^{\otimes k}$ lies in
$(E^{\sigma})^{\otimes k}$ and that $L_{\eta^{\otimes k}}^{\ast}$ maps
$(E^{\sigma})^{\otimes k}\otimes H$ into $H$ in the obvious way (and, for
$k=0$, this is $\iota_{H}$). Then, for every $X\in H^{\infty}(E)$,
\[
X(\eta^{\ast})=\sum_{k=0}^{\infty}L_{\eta^{\otimes k}}^{\ast}P_{k}U^{\ast
}(X\otimes I)U\iota_{H}%
\]
where $U$ is as defined in Lemma \ref{tensorsum}.
\end{proposition}

\section{Nevanlinna-Pick Theorem}

Our goal in this section is to present a generalization of the Nevanlinna-Pick
theorem. First, recall the classical theorem.

\begin{theorem}
\label{Npclassical}. Let $z_{1},\ldots z_{k} \in\mathbb{C}$ with $|z_{i}|<1$
and $w_{1},\ldots w_{k} \in\mathbb{C}$. Then the following conditions are equivalent.

\begin{itemize}
\item[(1)] There is a function $f\in H^{\infty}(\mathbb{T})$ with $\norm{f}
\leq1$ such that $f(z_{i})=w_{i}$ for all $i$.

\item[(2)]
\[
\left(  \frac{1-w_{i}\overline{w_{j}}}{1-z_{i}\overline{z_{j}}} \right)  \geq0
.
\]
\end{itemize}
\end{theorem}

Since we are able to view elements of $H^{\infty}(E)$ as functions on the open
unit ball of $E^{\sigma}$, it makes sense to seek necessary and sufficient
conditions for finding an element $X\in H^{\infty}(E)$ with norm less or equal
$1$ whose values at some prescribed points, $\eta_{1},\ldots\eta_{k}$, in that
open unit ball are prescribed operators $C_{1},\ldots C_{k}$ in $B(H)$. To
state our conditions we need some notation. For operators $B_{1}$, $B_{2}$ in
$B(H)$ we write $Ad(B_{1},B_{2})$ for the map, from $B(H)$ to itself, mapping
$S$ to $B_{1}SB_{2}^{\ast}$. Also, for elements $\eta_{1},\eta_{2}$ in
$\mathbb{D}(E^{\sigma})$, we let $\theta_{\eta_{1},\eta_{2}}$ denote the map,
from $\sigma(M)^{\prime}$ to itself, that sends $a$ to $\langle\eta_{1}%
,a\eta_{2}\rangle$. Then our generalization of the Nevanlinna-Pick theorem may
be formulated as follows.

\begin{theorem}
\label{NP}Let $\sigma$ be a faithful normal representation of $M$ on $H$. Fix
$\;\eta_{1},\ldots\eta_{k}\in E^{\sigma}$ with $\norm{\eta_i} <1$ and
$B_{1},\ldots B_{k},C_{1},\ldots C_{k}\in B(H)$. Then the following conditions
are equivalent

\begin{itemize}
\item[(1)] There exists an $X\in H^{\infty}(E)$ with $\norm{X} \leq1$ such
that $B_{i}X(\eta_{i}^{\ast})=C_{i}$ for all $i$.

\item[(2)] The map from $M_{k}(\sigma(M)^{\prime})$ into $M_{k}(B(H))$ defined
by the $k\times k$ matrix
\[
\left(  (Ad(B_{i},B_{j})-Ad(C_{i},C_{j}))\circ(id-\theta_{\eta_{i},\eta_{j}%
})^{-1} \right)
\]
is completely positive.
\end{itemize}
\end{theorem}

\begin{remark}
\label{choi} If $M=B(H)$ (and, then $\sigma(M)^{\prime}=\mathbb{C}I$),
condition (2) becomes
\[
\left(  \frac{B_{i}B_{j}^{*}-C_{i}C_{j}^{*}}{1-\langle\eta_{i},\eta_{j}%
\rangle} \right)  \geq0.
\]
This follows easily from a result of M. D. Choi (\cite{C75}).
\end{remark}

For the complete proof of Theorem~\ref{NP} we refer the reader to
\cite[Theorem 5.3]{MSNP}. Here we just remark that in order to prove that (1)
implies (2) one uses the complete positivity condition of (2) to construct a
subspace $\mathcal{M}\subseteq\mathcal{F}(E^{\sigma})\otimes H$ that is
invariant under $Ind(\iota)(H^{\infty}(E^{\sigma}))^{\ast}$ and a contraction
$R$ that commutes with the restriction of $Ind(\iota)(H^{\infty}(E^{\sigma
}))^{\ast}$to $\mathcal{M}$. Then it is possible to apply the commutant
lifting theorem of \cite[Theorem 4.4]{MS98} to $R^{\ast}$ to get a contraction
on $\mathcal{F}(E^{\sigma})\otimes H$ that commutes with $Ind(\iota
)(H^{\infty}(E^{\sigma}))$. An application of Theorem~\ref{commutant}
completes the proof.

The following is a consequence of Theorem \ref{NP}. It may be viewed as a
generalization of the classical Schwartz's lemma.

\begin{theorem}
\label{Theorem3.5}Suppose an element $X$ of $H^{\infty}(E)$ has norm at most
one and satisfies the equation $X(0)=0$. Then for every $\eta^{\ast}%
\in\mathbb{D}((E^{\sigma})^{\ast})$ the following assertions are valid:

\begin{enumerate}
\item If $a$ is a nonnegative element in $\sigma(M)^{\prime}$, and if
$\langle\eta,a\cdot\eta\rangle\leq a$, then
\[
X(\eta^{\ast})aX(\eta^{\ast})^{\ast}\leq\langle\eta,a\cdot\eta\rangle.
\]

\item If $\eta^{\otimes k}$ denotes the element $\eta\otimes\eta\otimes
\cdots\otimes\eta\in E^{\otimes k}$, then
\[
X(\eta^{\ast})\langle\eta^{\otimes k},\eta^{\otimes k}\rangle X(\eta^{\ast
})^{\ast}\leq\langle\eta^{\otimes k+1},\eta^{\otimes k+1}\rangle.
\]

\item $X(\eta^{\ast})X(\eta^{\ast})^{\ast}\leq\langle\eta,\eta\rangle$.
\end{enumerate}
\end{theorem}

We now illustrate how to apply Theorem~\ref{NP} in various settings.

\begin{example}
When $M=H=E=\mathbb{C}$, we obtain Theorem~\ref{Npclassical}.
\end{example}

\begin{example}
If  $M=E=\mathbb{C}$ and if $H$ is arbitrary, then $E^{\sigma}=B(H)$ and
Theorem \ref{NP} yields the following result.

\begin{theorem}
Given $T_{1},\ldots T_{k} \in B(H), \norm{T_i} <1$ and $B_{1},\ldots
,B_{k},C_{1}, \ldots C_{k}$ in $B(H)$. then the following conditions are equivalent.

\begin{itemize}
\item[(1)] There exists a function $f \in H^{\infty}(\mathbb{T})$ with
$\norm{f} \leq1$ and $B_{i}f(T_{i})=C_{i}$.

\item[(2)] The map defined by the matrix $(\phi_{ij})$ is completely positive
where
\[
\phi_{ij}(A)=\sum_{k=0}^{\infty} (B_{i}T_{i}^{k}AT_{j}^{*k}B_{j}^{*}%
-C_{i}T_{i}^{k}AT_{j}^{*k}B_{j}) .
\]
\end{itemize}
\end{theorem}
\end{example}

\begin{example}
Assume $M=B(H)=E$. Then $M^{\prime}=\mathbb{C}I$ and $E^{\sigma}=\mathbb{C}$
and Theorem \ref{NP} specializes to the following.

\begin{theorem}
Given $z_{1},\ldots z_{k} \in\mathbb{D}$ and $B_{1},\ldots,B_{k},C_{1},\ldots
C_{k}$ in $B(H)$, then the following conditions are equivalent.

\begin{itemize}
\item[(1)] There exists $G\in H^{\infty}(\mathbb{T})\otimes B(H)$ with
$\norm{G} \leq1$ such that $B_{i}G(z_{i})=C_{i}$ for all $i$.

\item[(2)]
\[
\left(  \frac{B_{i}B_{j}^{*}-C_{i}C_{j}^{*}}{1-z_{i}\overline{z_{j}}} \right)
\geq0 .
\]
\end{itemize}
\end{theorem}
\end{example}

\begin{example}
Set $M=B(H)$ and $E=C_{n}(B(H))$ (that is, $E$ is a column of $n$ copies of
$B(H)$). Then $M^{\prime}=\mathbb{C}I$, $E^{\sigma}=\mathbb{C}^{n}$ and
Theorem \ref{NP} yields the following theorem due to Davidson and Pitts
\cite{DP}, Arias and Popescu \cite{AP00} and Popescu \cite{Po01}.

\begin{theorem}
Given $\eta_{1},\ldots\eta_{k}$ in the open unit ball of $\mathbb{C}^{n}$ and
$C_{1},\ldots C_{k}\in B(H)$, then the following conditions are equivalent.

\begin{itemize}
\item[(1)] There is a $Y\in B(H)\otimes\mathcal{L}_{n}$ with $\norm{Y} \leq1$
such that $(\eta_{i}^{\ast}\times id)(Y)=C_{i}$ for all $i$.

\item[(2)]
\[
\left(  \frac{I-C_{i}C_{j}^{\ast}}{1-\langle\eta_{i},\eta_{j}\rangle}\right)
\geq0.
\]
Moreover, if, for all $i$, the $C_{i}$ all lie in some von Neumann subalgebra
$N\subseteq B(H)$, then $Y$ can be chosen in $N\otimes\mathcal{L}_{n}$.
\end{itemize}
\end{theorem}
\end{example}

Our final example of this section concerns interpolation for nest algebras.
The first interpolation result for nest algebras was proved by Lance
(\cite{L69}). It was later generalized by Anoussis (\cite{An92}) and by
Katsoulis, Moore and Trent (\cite{KMT93}). A related result was proved by Ball
and Gohberg (\cite{BG85}). The results we present below recapture the results
of \cite{KMT93}.

\begin{theorem}
Let $\mathcal{N}$ be a nest of projections in $B(H)$ and fix $B,C$ in $B(H)$.
Then the following conditions are equivalent.

\begin{itemize}
\item[(1)] There exists an $X\in Alg\mathcal{N}$ with $\norm{X} \leq1$ and
$BX=C$.

\item[(2)] For all projections $N \in\mathcal{N}$, $CNC^{*}\leq BNB^{*}$.
\end{itemize}
\end{theorem}

The ``vector version'' of this theorem is the following.

\begin{corollary}
Let $\mathcal{N}$ be a nest in $B(H)$ and fix $\;u_{1},\ldots u_{k}%
,v_{1},\ldots v_{k}$ in $H$. Then the following conditions are equivalent.

\begin{itemize}
\item[(1)] There exists $X \in Alg\mathcal{N}$ with $\norm{X} \leq1$ and
$Xu_{i}=v_{i}$ for all $i$.

\item[(2)] For all $N\in\mathcal{N}$,
\[
\left(  \langle N^{\perp}v_{i},N^{\perp}v_{j}\rangle\right)  \leq\left(
\langle N^{\perp}u_{i},N^{\perp}u_{j}\rangle\right)  \text{,}%
\]
where $N^{\perp}$ denotes $I-N$.
\end{itemize}
\end{corollary}

These results are not immediate corollaries of Theorem \ref{NP} because, for a
general nest $\mathcal{N}$, $Alg\mathcal{N}$ is not of the form $H^{\infty
}(E)$. However, when $\mathcal{N}$ is finite, $Alg\mathcal{N}$ is a Hardy
Algebra by Example \ref{ex0.4}. In this case, the conclusions are fairly
straight forward computations. The case of general nests is then handled by
approximation techniques along the lines of \cite{L69} and \cite{Ar75}. Full
details may be found in \cite[Theorem 6.8 and Corollary 6.9]{MSNP}.

\section{Schur class operator functions and realization}

In this section we relate the complete positivity condition of
Theorem~\ref{NP} to the concept of a Schur class function. As mentioned in the
introduction, this may be viewed as part of a general program to find
equivalences between canonical model theory, ``non commutative'' systems
theory and scattering theory. The results below are proved in \cite{MSSchur}.

We start with the following definition.

\begin{definition}
\label{Cpker} Let $S$ be a set, $A$ and $B$ be two $C^{\ast}$-algebras and
write $\mathcal{B}(A,B)$ for the space of bounded linear maps from $A$ to $B$.
A function
\[
K:S\times S\rightarrow\mathcal{B}(A,B)
\]
will be called a \emph{completely positive definite kernel} (or a
\emph{CPD-kernel} ) if, for all choices of $s_{1},\ldots,s_{k}$ in $S$, the
map
\[
K^{(k)}:(a_{ij})\mapsto(K^{(k)}(s_{i},s_{j})(a_{ij}))
\]
from $M_{k}(A)$ to $M_{k}(B)$ is completely positive.
\end{definition}

This concept of CPD-kernels was studied in \cite{BBLS} (see, in particular,
Lemma 3.2.1 there for conditions on $K$ that are equivalent to being a CPD-kernel).

\begin{definition}
\label{schur} Let $\mathcal{E}$ be a Hilbert space and $Z: \mathbb{D}%
((E^{\sigma})^{*}) \rightarrow B(\mathcal{E})$ be a $B(\mathcal{E})$-valued
function. Then $Z$ is said to be a \emph{Schur class operator function} if
\[
K(\eta^{*},\zeta^{*})=(id-Ad(Z(\eta^{*}),Z(\zeta^{*}))\circ(id-\theta
_{\eta,\zeta})^{-1}%
\]
is a CPD-kernel on $\mathbb{D}((E^{\sigma})^{*})$. (We use here the notation
set for Theorem~\ref{NP}).
\end{definition}

Note that, when $M=E=B(\mathcal{E})$ and $\sigma$ is the identity
representation of $B(\mathcal{E})$ on $\mathcal{E}$, $\sigma(M)^{\prime}$ is
$\mathbb{C}I_{\mathcal{E}}$, $E^{\sigma}$ is isomorphic to $\mathbb{C}$ and
$\mathbb{D}((E^{\sigma})^{\ast})$ may be identified with the open unit ball
$\mathbb{D}$ of $\mathbb{C}$. In this case the definition above recovers the
classical Schur class functions. More precisely, these functions are usually
defined as analytic functions $Z$ from $\mathbb{D}$ into the closed unit ball
of $B(\mathcal{E})$ but it is known that this is equivalent to the positivity
of the Pick kernel $k_{Z}(z,w)=(I-Z(z)Z(w)^{\ast})(1-z\bar{w})^{-1}$. The
argument of \cite[Remark 5.4]{MSNP} shows that the positivity of this kernel
is equivalent, in this case, to the condition of Definition~\ref{schur}.

Note that it follows from Theorem~\ref{NP} that every operator in the closed
unit ball of $H^{\infty}(E)$ determines (by point evaluation) a Schur class
operator function. In fact we have the following result whose proof may be
found in \cite{MSSchur}.

\begin{theorem}
\label{schurealiz} (\cite{MSSchur}) Let $M$ be a von Neumann algebra, let $E$
be a $W^{\ast}$-correspondence over $M$ and let $\sigma$ be a faithful normal
representation of $M$ on a Hilbert space $\mathcal{E}$. For a function
$Z:\mathbb{D}((E^{\sigma})^{\ast})\rightarrow B(\mathcal{E})$, the following
conditions are equivalent.

\begin{enumerate}
\item[(1)] $Z$ is a Schur class operator function.

\item[(2)] There is an $X$ in the closed unit ball of $H^{\infty}(E)$ such
that $X(\eta^{\ast})=Z(\eta^{\ast})$ for all $\eta\in\mathbb{D}(E^{\sigma})$.

\item[(3)] (Realization) There is a Hilbert space $H$, a normal representation
$\tau$ of $N :=\sigma(M)^{\prime}$ on $H$ and operators $A,B,C $ and $D$ such that

\begin{enumerate}
\item[(i)] $A\in B(\mathcal{E})$,$B\in B(H,\mathcal{E})$, $C\in B(\mathcal{E}%
,H)$ and $D\in B(H, E^{\sigma}\otimes H)$.

\item[(ii)] $A,B,C$ and $D$ intertwine the actions of $N$ (on the relevant spaces).

\item[(iii)] The operator
\[
V :=\left(
\begin{array}
[c]{cc}%
A & B\\
C & D
\end{array}
\right)  : \left(
\begin{array}
[c]{c}%
\mathcal{E}\\
H
\end{array}
\right)  \rightarrow\left(
\begin{array}
[c]{c}%
\mathcal{E}\\
E^{\sigma}\otimes H
\end{array}
\right)
\]
is a coisometry.

\item[(iv)] For every $\eta\in\mathbb{D}(E^{\sigma})$,
\[
Z(\eta^{*})=A+B(I-L_{\eta}^{*}D)^{-1}L_{\eta}^{*}C
\]
where $L_{\eta}:H \rightarrow E^{\sigma}\otimes H$ is defined by $L_{\eta
}h=\eta\otimes h$.
\end{enumerate}
\end{enumerate}
\end{theorem}

Note that $X$ in part (2) of the Theorem is not necessarily unique. (Although,
as shown in \cite{MSSchur}, it is possible to choose $\sigma$ such that the
choice of $X$ will be unique).

One may apply the techniques developed (in \cite{MSSchur}) for the proof of
the Theorem \ref{schurealiz} to establish the following extension result.

\begin{proposition}
\label{extension} Every function defined on a subset $\Omega$ of the open unit
ball $\mathbb{D}((E^{\sigma})^{\ast})$ with values in some $B(\mathcal{E})$
such that the associated kernel (defined on $\Omega\times\Omega$) is a
CPD-kernel may be extended to a Schur class operator function (defined on all
of $\mathbb{D}((E^{\sigma})^{\ast})$).
\end{proposition}


\begin{thebibliography}{9}                                                                                                %
\bibitem {AM}J. Agler and J. McCarthy, \emph{Pick interpolation and Hilbert
function spaces}. Graduate Studies in Mathematics, vol. 44. Amer. Math. Soc.,
Providence (2002).

\bibitem {ADD90}D. Alpay, P. Dewilde and H. Dym, \emph{Lossless inverse
scattering and reproducing kernels for upper triangular operators}, Operator
Theory: Adv. Appl., Birkhauser Verlag, Basel 47 (1990), 61-133.

\bibitem {ADP}D. Alpay, A Dijksma and Y. Peretz, \emph{Nonstationary analogues
of the Herglotz representation theorem: the discrete case}, J. Funct. Anal. 66
(1999), 85-129.

\bibitem {AD90}C. Anantharaman-Delaroche, \emph{On completely positive maps
defined by an irreducible correspondence}, Canad. Math. Bull. 33 (1990), 434-441.

\bibitem {An92}M. Anoussis, \emph{Interpolating operators in nest algebras,
}Proc. Amer. Math. Soc. 114 (1992), 707--710.

\bibitem {AP00}A. Arias and G. Popescu, \emph{Noncommutative interpolation and
Poison transforms, }Israel J. Math. 115 (2000), 205--234.

\bibitem {Ar75}Wm. Arveson, \emph{Interpolation problems in nest algebras, }J.
Funct. Anal. \textbf{3 }(1975), 208--233.

\bibitem {Ar89}W.B. Arveson, \emph{Continuous analogues of Fock space}, Mem.
Amer. Math. Soc. 80 (1989).

\bibitem {BDH88}M. Baillet, Y. Denizeau and J.-F. Havet, \emph{Indice d'une
esperance conditionelle}, Comp. Math. 66 (1988), 199-236.

\bibitem {B}J. Ball, \emph{Linear systems, operator model theory and
scattering: multivariable generalizations}, in Operator Theory and its
applications (Winnipeg, 1998) (Ed. A.G. Ramm, P.N. Shivakumar and A.V.
Strauss), Fields Inst. Comm. vol. 25, Amer. Math. Soc., Providence, 2000, 151-178.

\bibitem {BG85}J. Ball and I. Gohberg, \emph{A commutant lifting theorem for
triangular matrices with diverse applications}, Integral Equat. Operator
Theory 8 (1985), 205--267.

\bibitem {BTV}J. Ball, T. Trent and V. Vinnikov, \emph{Interpolation and
commutant lifting for multipliers on reproducing kernel Hilbert spaces}. Preprint.

\bibitem {BV}J. Ball and V. Vinnikov, \emph{Functional models for
representations of the Cuntz algebra}. Preprint.

\bibitem {BBLS}S.D. Barreto, B.V.R. Bhat, V. Liebscher and M. Skeide,
\emph{Type I product systems of Hilbert modules}, to appear in J. Functional Anal.


\bibitem {hB88}H. Bercovici, \emph{Operator theory and arithmetic in
}$H^{\infty}$. Mathematical Surveys and Monographs, 26. American Mathematical
Society, Providence, RI, 1988.

\bibitem {C75}M.D. Choi, \emph{Completely positive linear maps on complex
matrices}, Lin. Alg. Appl. 10 (1975), 285-290.

\bibitem {DP}K. Davidson and D. Pitts, \emph{The algebraic structure of
non-commutative analytic Toeplitz algebras}, Math. Ann. 311 (1998), 275-303.

\bibitem {pG72}P. Gabriel, \emph{Unzerlegbare Darstellungen I, }Manuscr. Math.
\textbf{6 }(1972), 71--103.

\bibitem {GR91}P. Gabriel and A. V. Roiter, \emph{Representations of
Finite-Dimensional Algebras, }Algebra VIII, Encyclopaedia of Mathematical
Sciences, Vol. 73, Springer-Verlag, 1991.

\bibitem {gH47}G. Hochschild, \emph{On the structure of algebras with nonzero
radical, }Bull. Amer. Math. Soc. \textbf{53 }(1947), 369--377.

\bibitem {KMT93}E.G.\ Katsoulis, R.L. Moore and T.T. Trent,
\emph{Interpolation in nest algebras and applications to operator corona
theorems, }J. Operator Th. 29 (1993), 115--123.

\bibitem {KP}D. Kribs and S. Power, \emph{Free semigroupoid algebras}, Preprint.

\bibitem {L69}E. C. Lance, \emph{Some properties of nest algebras, }Proc.
London Math. Soc. (3) \textbf{19 }(1969), 45--68.

\bibitem {L94}E.C. Lance , \emph{Hilbert $C^{\ast}$-modules, A toolkit for
operator algebraists}, London Math. Soc. Lecture Notes series 210 (1995).
Cambridge Univ. Press.

\bibitem {MM}M. McAsey and P.S. Muhly, \emph{Representations of
non-self-adjoint crossed products}, Proc. London Math. Soc. 47 (1983), 128-144.

\bibitem {jM}J. Mingo, \emph{The correspondence associated to an inner
completely positive map}, Math. Ann. 284 (1989), 121-135.

\bibitem {pM97}P.S. Muhly, \emph{A finite-dimensional introduction to operator
algebra} in \emph{Operator algebras and applications (Samos, 1996)}, 313--354,
NATO Adv. Sci. Inst. Ser. C Math. Phys. Sci., 495, Kluwer Acad. Publ.,
Dordrecht, 1997.

\bibitem {MS98}P.S. Muhly and B. Solel, \emph{Tensor algebras over $C^{\ast}%
$-correspondences (Representations, dilations and $C^{\ast}$-envelopes)}, J.
Funct. Anal. 158 (1998), 389-457.

\bibitem {MS99}P.S. Muhly and B. Solel , \emph{Tensor algebras, induced
representations, and the Wold decomposition}, Canad. J. Math. 51 (1999), 850-880.

\bibitem {MSQMP}P.S. Muhly and B. Solel, \emph{Quantum Markov processes
(correspondences and dilations)}, Int. J. Math. 13 (2002), 863-906.

\bibitem {MSC}P.S. Muhly and B. Solel, \emph{The curvature and index of
completely positive maps}, Proc. London Math. Soc. 87 (2003), 748-778.

\bibitem {MSNP}P.S. Muhly and B. Solel, \emph{Hardy algebras, $W^{\ast}%
$-correspondences and interpolation theory}, to appear in Math. Ann.

\bibitem {MSCM}P.S. Muhly and B. Solel, \emph{On canonical models for
representations of Hardy algebras}. In preparation.

\bibitem {MSSchur}P.S. Muhly and B. Solel, \emph{Schur class operator
functions associated with a $W^{*}$-correspondence}. In preparation.

\bibitem {MSS03}P.S. Muhly, M. Skeide and B. Solel, \emph{Representations of
$\mathcal{B}^{a}(E)$, commutants of von Neumann bimodules, and product systems
of Hilbert modules}. In preparation.

\bibitem {NF66}B. Sz-Nagy and C. Foias, \emph{Analyse Harmonique des
Operateurs de L'espace de Hilbert}, Akademiai Kiado (1966).

\bibitem {wP73}W. Paschke, \emph{Inner product modules over }$B^{\ast}%
$\emph{-algebras}, Trans. Amer. Math. Soc. \textbf{182 }(1973), 443--468.

\bibitem {Pe}J. Peters, \emph{Semi-crossed products of C*-algebras}, J. Funct.
Anal. 59 (1984), 498-534.

\bibitem {Pi}M. Pimsner, \emph{A class of $C^{\ast}$-algebras generalyzing
both Cuntz-Krieger algebras and crossed products by $\mathbb{Z}$}, in
\textit{Free Probability Theory}, D. Voiculescu, Ed., Fields Institute Comm.
12, 189-212, Amer. Math. Soc., Providence, 1997.

\bibitem {P86}S. Popa, \emph{Correspondences}, Preprint (1986).

\bibitem {Po89}G. Popescu, \emph{Characteristic functions for infinite
sequences of noncommuting operators}, J. Oper. Theory 22 (1989), 51-71.

\bibitem {Po91}G. Popescu, \emph{von Neumann inequality for $B(\mathcal{H}%
^{n})_{1}$}, Math. Scand. 68 (1991), 292-304.

\bibitem {Po95}G. Popescu, \emph{Functional calculus for noncommuting
operators, }Mich. Math. J. \textbf{42 }(1995), 345--356.

\bibitem {Po01}G. Popescu, \emph{Commutant lifting, tensor algebras and
functional calculus, }Proc. Edinburg Math. Soc. 44 (2001), 389--406.

\bibitem {Rie}M.A. Rieffel, \emph{Induced representations of $C^{\ast}%
$-algebras}, Adv. in Math. 13 (1974), 176-257.

\bibitem {S67}D. Sarason, \emph{Generalized interpolation in $H^{\infty}$},
Trans. Amer. Math. Soc. 127 (1967), 179-203.


\end{thebibliography}
\end{document}